\documentclass[12pt]{article}
\usepackage{setspace}
\usepackage{amssymb}
\usepackage{amsmath}
\usepackage{gdgspace}

\newcommand{\reff}[1]{\ref{#1}}

\numberwithin{equation}{section}

\begin{document}

\vskip .5cm

\begin{center}
\textbf{A NOTE ON $q$-ANALOGUE OF S\'{A}NDOR'S FUNCTIONS}
\end{center}
%\vskip .5cm
\begin{center}
by
\end{center}
\begin{center}
\textbf{ Taekyun Kim$^{1}$,  C. Adiga$^{2}$   and Jung Hun Han$^{2}$}\\
\vskip.5cm
$^{1}$Department of Mathematics Education\\
Kongju National University, Kongju 314-701\\
South Korea\\
$\text{e-mail:}$$ \text{tkim@kongju.ac.kr /
tkim64@hanmail.net}$\\\vskip .5cm
$^{2}$Department of Studies in Mathematics\\
University of Mysore, Manasagangotri\\
Mysore 570006, India\\
$\text{e-mail:c-adiga@hotmail.com}$ \vskip .5cm
\end{center}
\vskip .5cm

\begin{center}
Dedicated to Sun-Yi Park on 90th birthday
\end{center}
\vskip .5cm

%------------------------------------------------------------------

\begin{center}
\textbf{ABSTRACT}
\end{center}
\begin{spacing}{1.5}

%------------------------------------------------------------------------------------------

The additive analogues of Pseudo-Smarandache, Smarandache-simple
functions and their duals have been recently studied by J.
S\'{a}ndor. In this note, we obtain $q$-analogues of S\'{a}ndor's
theorems \cite{s-2}.

\begin{tabbing}
Smarandache-simple function,
Asymtotic formula.
\end{tabbing}
\begin{tabbing}
\end{tabbing}

%------------------------------------------------------------------------------------------

\section{Introduction}
\label{intro} \setcounter{equation}{0} \quad \quad The additive
analogues of Smarandache functions $S$ and $S_*$ have been
introduced by S\'{a}ndor
\cite{s-1} as follows:\\
$$S(x)=\text{min}\{m\in N: x\leq m!\},\quad x\in (1,\infty),$$
and
$$S_*(x)=\text{max}\{m\in N: m!\leq x\},\quad x\in [1,\infty).$$
He has studied many important properties of $S_*$ relating to
continuity, differentiability and Riemann integrability and also
proved the following theorems:

{\bf Theorem 1.1}
$$S_*(x) \sim \frac{\text{log} x}{\text{log} \text{log} x} \quad
(x\rightarrow \infty).$$

{\bf Theorem 1.2} The series
$$\sum_{n=1}^{\infty}\frac{1}{n(S_*(n))^\alpha}$$
is convergent for $\alpha >1$ and divergent for $\alpha\leq 1$.\\

In \cite{a-k}, Adiga and Kim have obtained generalizations of
Theorems 1.1 and 1.2 by the use of Euler's gamma function.
Recently Adiga-Kim-Somashekara-Fathima  \cite{a-k-s-f} have
established a $q$-analogues of these results on employing the
$q$-analogue of Stirling's formula. In \cite{s-2}, S\'{a}ndor
defined the additive analogues of Pseudo-Smarandache,
Smarandache-simple functions and their duals as follows:
$$Z(x)=\text{min}\left\{m\in N: x\leq \frac{m(m+1)}{2}\right\},\quad x\in
(0,\infty),$$
$$Z_*(x)=\text{max}\left\{m\in N: \frac{m(m+1)}{2}\leq x\right\},\quad x\in
[1,\infty),$$
$$P(x)=\text{min}\{m\in N: p^x\leq m!\},\quad p>1, x\in (0,\infty),$$
and
$$P_*(x)=\text{max}\{m\in N: m!\leq p^x\},\quad p>1, x\in [1,\infty).$$
He has also proved the following theorems:

{\bf Theorem 1.3}
$$Z_*(x) \sim \frac{1}{2}\sqrt{8x+1} \quad (x\rightarrow \infty).$$

{\bf Theorem 1.4} The series
$$\sum_{n=1}^{\infty}\frac{1}{(Z_*(n))^\alpha}$$
is convergent for $\alpha >2$ and divergent for $\alpha\leq 2$.
The series
$$\sum_{n=1}^{\infty}\frac{1}{n(Z_*(n))^\alpha}$$
is convergent for all $\alpha >0$.

{\bf Theorem 1.5}
$$\text{log}P_*(x) \sim \text{log}x \quad (x\rightarrow \infty).$$

{\bf Theorem 1.6} The series
$$\sum_{n=1}^{\infty}\frac{1}{n} \left(\frac{ \text{log} \text{log} n}
{\text{log} P_*(n)}\right)^\alpha$$
is convergent for $\alpha >1$ and divergent for $\alpha\leq 1$.\\

The main purpose of this note is to obtain $q$-analogues of
S\'{a}ndor's Theorems 1.3 and 1.5. In what follows, we make use of
the following notations and definitions. F. H. Jackson defined a
$q$-analogue of the gamma function which extends the $q$-factorial
$$(n!)_q = 1(1+q)(1+q+q^2)\cdots (1+q+q^2+\cdots +q^{n-1}),\quad\text{cf
\cite{k}},$$
which becomes the ordinary factorial as $q\rightarrow
1$. He defined the $q$-analogue of the gamma function as
$$\Gamma_q(x)= \frac{(q;q)_\infty}{(q^x;q)_\infty}(1-q)^{1-x},\quad 0<q<1,$$

and

$$\Gamma_q(x)= \frac{(q^{-1};q^{-1})_\infty}{(q^{-x};q^{-1})_\infty}
(q-1)^{1-x}q^{x \choose 2},\quad q>1,$$ where
$$(a;q)_\infty = \prod^{\infty}_{n=0}(1-aq^n).$$

It is well known that $\Gamma_q(x)\rightarrow \Gamma(x)$ as
$q\rightarrow 1$, where $\Gamma(x)$ is the ordinary gamma
function.

%--------------------------------------------------------------------------------------------------

\section{Main Theorems}

We now define the $q$-analogues of $Z$ and $Z_*$ as follows:
$$Z_q(x)=\text{min}\left\{\frac{1-q^m}{1-q}: x\leq
\frac{\Gamma_q(m+2)}{2\Gamma_q(m)}\right\}
, \quad m \in N, x\in (0,\infty),$$ and
$$Z_q^*(x)=\text{max}\left\{\frac{1-q^m}{1-q}:
\frac{\Gamma_q(m+2)}{2\Gamma_q(m)}\leq x\right\}
, \quad m \in N, x\in
\left[\frac{\Gamma_q(3)}{2\Gamma_q(1)},\infty\right),$$ where
$0<q<1$. Clearly, $Z_q(x)\rightarrow Z(x)$ and
$Z_q^*(x)\rightarrow Z_*(x)$ as $q\rightarrow 1^{-}$. From the
definitions of $Z_q$ and $Z_q^*$, it is clear that

\begin{equation}
Z_q(x) = \left\{\begin{array}{ll} 1, \quad \mbox{if \quad$ \quad x
\in
\left(0,\frac{\Gamma_q(3)}{2\Gamma_q(1)}\right]$} \\
\frac{1-q^{m}}{1-q} , \quad \mbox{if} \quad x \in \left(
\frac{\Gamma_q(m+1)}{2\Gamma_q(m-1)},
\frac{\Gamma_q(m+2)}{2\Gamma_q(m)}\right], m \geq 2,\end{array}
\label{zq-1} \right]
\end{equation}

and

\begin{equation}
Z_q^*(x) =  \frac{1-q^m}{1-q} \quad \text{if} \quad x \in
\left[\frac{\Gamma_q(m+2)}{2\Gamma_q(m)},\frac{\Gamma_q(m+3)}{2\Gamma_q(m+1)}\right).
\label{zq*-1}
\end{equation}

Since
$$\frac{1-q^{m-1}}{1-q} \leq \frac{1-q^m}{1-q} =
\frac{1-q^{m-1}}{1-q} + q^{m-1} \leq \frac{1-q^{m-1}}{1-q} +1,$$

(\reff{zq-1}) and (\reff{zq*-1}) imply that for $x\geq
\frac{\Gamma_q(3)}{2\Gamma_q(1)}$,
$$Z_q^*(x) \leq Z_q(x) \leq Z_q^*(x)+1 .$$
Hence it suffices to study the function $Z_q^*$. We now prove our
main theorems.

%----------------------------------------------------------------------------------

{\bf Theorem 2.1} If $0<q<1$, then
$$\frac{\sqrt{1+8xq}-(1+2q)}{2q^2}< Z_q^*(x) \leq \frac{\sqrt{1+8xq}-1}{2q}
, \quad x\geq \frac{\Gamma_q(3)}{2\Gamma_q(1)}.$$

{\bf Proof.} If
\begin{equation}
\frac{\Gamma_q(k+2)}{2\Gamma_q(k)} \leq x <
\frac{\Gamma_q(k+3)}{2\Gamma_q(k+1)}, \label{gamma-1}
\end{equation}

then
$$Z_q^*(x)=\frac{1-q^k}{1-q}$$
and
\begin{equation}
(1-q^k)(1-q^{k+1})-2x(1-q)^2 \leq 0<
(1-q^{k+1})(1-q^{k+2})-2x(1-q)^2. \label{gamma-2}
\end{equation}

Consider the functions $f$ and $g$ defined by
$$f(y)= (1-y)(1-yq)- 2x(1-q)^2$$
and
$$g(y)=(1-yq)(1-yq^2)-2x(1-q)^2.$$
Note that $f$ is monotonically decreasing for
$y\leq\frac{1+q}{2q}$ and $g$ is strictly decreasing for $y <
\frac{1+q}{2q^2}$. Also $f(y_1)=0=g(y_2)$ where
$$y_1=\frac{(1+q)-(1-q)\sqrt{1+8xq}}{2q},$$
$$y_2=\frac{(q+q^2)-q(1-q)\sqrt{1+8xq}}{2q^3}.$$

Since $y_1< \frac{1+q}{2q}, y_2< \frac{1+q}{2q^2}$ and
$q^k<\frac{1+q}{2q}< \frac{1+q}{2q^2},$ from (\reff{gamma-2}), it
follows that
$$f(q^k) \leq f(y_1)=0 = g(y_2)< g(q^k).$$
Thus $y_1 \leq q^k < y_2$ and hence
$$\frac{1-y_2}{1-q} < \frac{1-q^k}{1-q} \leq \frac{1-y_1}{1-q}.$$
i.e.
$$\frac{\sqrt{1+8xq}-(1+2q)}{2q^2}< Z_q^*(x) \leq
\frac{\sqrt{1+8xq}-1}{2q}.$$ This completes the proof.

{\bf Remark.} Letting $q\rightarrow 1^-$ in the above theorem, we
obtain S\'{a}ndor's Theorem 1.3.\\
We define the $q$-analogues of $P$ and $P_*$ as follows:

$$P_q(x)=\text{min}\{m\in N: p^x\leq \Gamma_q(m+1)\}
, \quad p>1, x\in (0,\infty),$$ and

$$P_q^*(x)=\text{max}\{m\in N: \Gamma_q(m+1)\leq p^x\}
, \quad p>1, x\in [1,\infty),$$ where $0<q<1$. Clearly,
$P_q(x)\rightarrow P(x)$ and $P_q^*(x)\rightarrow P_*(x)$ as
$q\rightarrow 1^{-}$. From the definitions of $P_q$and $P_q^*$, we
have
$$P_q^*(x) \leq P_q(x) \leq P_q^*(x) +1.$$
Hence it is enough to study the function $P_q^*$.

%--------------------------------------------------------------------------------------

{\bf Theorem 2.2} If $0<q<1$, then
$$P_*(x)\sim \frac{x \text{log} p}{\text{log} \left(\frac{1}{1-q}\right)}
\quad (x\rightarrow \infty).$$

{\bf Proof.} If $\Gamma_q(n+1)\leq p^x < \Gamma_q(n+2)$, then
$$P_q^*(x)=n$$
and
\begin{equation}
\text{log} \Gamma_q(n+1)\leq \text{log} p^x < \text{log}
\Gamma_q(n+2). \label{gamma-3}
\end{equation}

But by the $q$-analogue of Stirling's formula established by Moak
\cite{m}, we have
\begin{equation}
\text{log} \Gamma_q(n+1) \sim \left(n+\frac{1}{2}\right)
\text{log} \left(\frac{q^{n+1}-1}{q-1}\right) \sim n \text{log}
\left(\frac{1}{1-q}\right). \label{gamma-4}
\end{equation}
Dividing (\reff{gamma-3}) throughout by $n \text{log}
(\frac{1}{1-q})$, we obtain
\begin{equation}
\frac{\text{log} \Gamma_q(n+1)}{n \text{log} (\frac{1}{1-q})} \leq
\frac{x \text{log} p}{P_q^*(x)\text{log} (\frac{1}{1-q})} <
\frac{\text{log} \Gamma_q(n+2)}{n \text{log}(\frac{1}{1-q})}  .
\label{gamma-5}
\end{equation}

Using (\reff{gamma-4}) in (\reff{gamma-5}), we deduce
$$\lim_{x\rightarrow \infty} \frac{x \text{log} p}{P_q^*(x) \text{log}
(\frac{1}{1-q})}=1.$$
This completes the proof.

%--------------------------------------------------------------------------------------

%{\bf Addresses:}\\

%Department of Studies in Mathematics\\
%University of Mysore, Manasagangotri, Mysore 570 006, India\\

\end{spacing}
\end{document}